\documentclass[11pt]{article} \textwidth=13cm \textheight=21cm
\usepackage[english]{babel}
\usepackage{pstricks,pst-node,pst-text,pst-3d,pst-plot}
\usepackage[active]{srcltx}

\begin{document}

\centerline{\Large{${\cal A}_{5}$ and the icosahedron rotations:}}
\centerline {\Large{seeing the facts, not only proving them.}}

\vspace{0.2cm}

\centerline{J. M. S. Sim\~{o}es-Pereira (*)}

\vspace{0.2cm}

\textit{Abstract: When compared with pure mathematicians, applied ones have a clear preference for proofs that go beyond a chain of reasonings and \emph{do exhibit the fact }to be proved. Here we exhibit the bijection between the 60 icosahedron rotations of the group $\cal R$ and the 60 permutations of the group ${\cal A}_{5}$.}

\vspace{0.2cm}

The motivation for this paper is \textit{educational}. Quite often, applied math students do not content themselves with the details of a proof; they want the details of \textit{the fact} that has been proved. This evokes the deep philosophical distinction between \emph{leibnizian} and \emph{cartesian} proofs as explained in Hacking  \cite{Hacking}. In a \emph{leibnizian} proof, we have a chain of steps, seen one by one; contrary to that, \emph{Descartes} claimed he wanted to have the entirety of a proof in his mind at once, an option which is clearly similar to our \emph{seeing the facts}.

A truly elementary example is the Pythagoras' theorem, the equality $c^{2}=a^{2}+b^{2}$ in a right-angled triangle where $c,\  a,\  b$ are the lengths of the hypotenuse and of the other two sides, respectively. Fig. 0 is used to exhibit this.

\vspace{0.15cm}

$$\psmatrix[mnode=circle,colsep=0.6cm,rowsep=0.6cm]1&2&&3&&1&2&&&3\\&&&4&&5&6&&&7\\8\\9&&0&y&&z&p&&&q
\endpsmatrix $$ \ncline{1,1}{1,2}\taput{a}\ncline{1,2}{1,4}\taput{b}\ncline{1,6}{1,7}\taput{a}\ncline{1,7}{1,10}\taput{b}
\ncline{1,1}{3,1}\tlput{b}\ncline{3,1}{4,1}\tlput{a}\ncline{1,4}{2,4}\trput{a}\ncline{2,4}{4,4}\trput{b}
\ncline{4,1}{4,3}\tbput{b}
\ncline{4,3}{4,4}\tbput{a}
\ncline{3,1}{1,2}\trput{c}\ncline{1,2}{2,4}\tlput{c}\ncline{2,4}{4,3}\tlput{c}\ncline{4,3}{3,1}\trput{c}
\ncline{2,6}{2,7}\ncline{2,7}{2,10}\ncline{4,6}{4,7}\tbput{a}\ncline{4,7}{4,10}\tbput{b}
\ncline{1,6}{2,6}\tlput{a}\ncline{2,6}{4,6}\tlput{b}\ncline{1,7}{2,7}\ncline{2,7}{4,7}\ncline{1,10}{2,10}\trput{a}
\ncline{2,10}{4,10}
\trput{b}
\ncline{2,7}{1,10}\ncline{2,6}{4,7}
\rput(3.9,2.6){$c^{2}$}\rput(9.7,2.1){$b^{2}$}\rput(8,3.8){$a^{2}$}

\vspace{0.15cm}

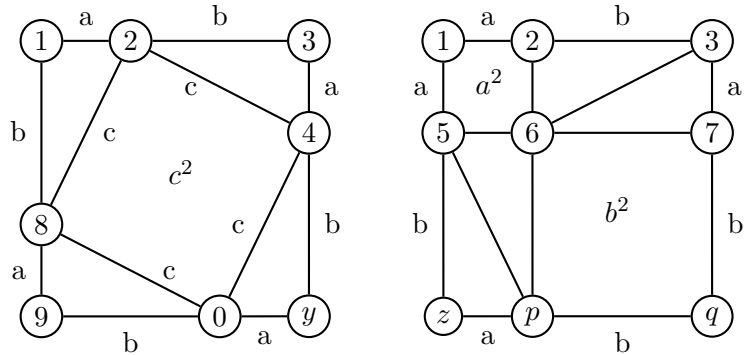
\centerline{Fig. 0: A graphic proof of the Pythagorean Theorem}

\vspace{0.15cm}

Looking at Fig. 0, it is clear that subtracting from the area of a square with side $a+b$, the total area of the four interior right-angled triangles whose sides are $a$, $b$ and $c$, we get $c^{2}$ on the left half of the picture, and we get $a^{2}+b^{2}$  on the right half of the picture. Obviously, this shows that  $c^{2}=a^{2}+b^{2}$.

In what follows, as another example, we take the isomorphism between the icosahedron rotations group $\cal R$ and the alternating group ${\cal A}_{5}$. In century old references like \cite{Burnside,Klein,MBD} or contemporary web sites like \cite{stack}, we find proofs of this result. Our purpose is not to give another one, although one is obtained here as a byproduct. Our purpose is twofold: \textit{to detail the bijection} between the elements of $\cal R$ and those of ${\cal A}_{5}$, that is to say, to give a complete list of this correspondence; and \emph{to exhibit how a set of generators} yields each one of the group elements. Applied math students have asked me to show them both things. They will also see that sometimes a lot of hard work is unavoidable!

\newpage $ $

\vspace{1.3cm}

\hspace{1.cm}\psset{xunit=0.09cm,yunit=0.09cm,arrowsize=3pt
2,arrowlength=3,shortput=tablr}
\cnodeput(0,0){222}{$2''$}\cnodeput(116,0){333}{$3''$}
\ncarc[arcangle=20]{222}{333} \cnodeput(58,-6){11}{$1'$}
\ncline{222}{11}\ncline{11}{333}\cnodeput(58,-22){1}{1}\ncline{11}{1}\cnodeput(58,-64){1m}{$1^{+}$}
\cnodeput(58,-100){111}{$1''$}\ncline{1m}{111}
\ncarc[arcangle=-20]{222}{111}\ncarc[arcangle=20]{333}{111}
\cnodeput(33,-18){2m}{$2^{+}$}\cnodeput(83,-18){3m}{$3^{+}$}
\ncline{11}{2m}\ncline{11}{3m}\ncline{1}{2m}\ncline{1}{3m}
\ncline{222}{2m}\ncline{333}{3m}
\cnodeput(49,-40){3}{3}\cnodeput(67,-40){2}{2}
\ncline{1}{2}\ncline{2}{3}\ncline{3}{1}
\ncline{1m}{3}\ncline{2}{1m}
\cnodeput(36,-48){33}{$3'$}\cnodeput(80,-48){22}{$2'$}
\ncline{33}{222}\ncline{33}{2m}\ncline{33}{3}\ncline{33}{1m}\ncline{33}{111}
\ncline{22}{333}\ncline{22}{3m}\ncline{22}{2}\ncline{22}{1m}\ncline{22}{111}
\ncline{3}{2m}\ncline{2}{3m} \rput(58,3){$\overrightarrow{T}$}
\rput(83,-9){$\overrightarrow{W}$}
\rput(33,-9){$\overrightarrow{V}$}
\rput(51,-15){$\overleftarrow{B}$}
\rput(65,-15){$\overleftarrow{Z}$}
\rput(49,-27){$\overleftarrow{D}$}
\rput(67,-27){$\overleftarrow{Y}$}
\rput(58,-33){$\overleftarrow{X}$}
\rput(58,-48){$\overleftarrow{T}$}
\rput(24,-22){$\overrightarrow{C}$}
\rput(92,-22){$\overrightarrow{A}$}
\rput(76,-36){$\overleftarrow{C}$}
\rput(40,-36){$\overleftarrow{A}$}
\rput(47,-49){$\overleftarrow{W}$}
\rput(69,-49){$\overleftarrow{V}$}
\rput(51,-71){$\overrightarrow{Z}$}
\rput(65,-71){$\overrightarrow{B}$}
\rput(28,-53){$\overrightarrow{Y}$}
\rput(88,-53){$\overrightarrow{D}$}

\vspace{9.7cm}

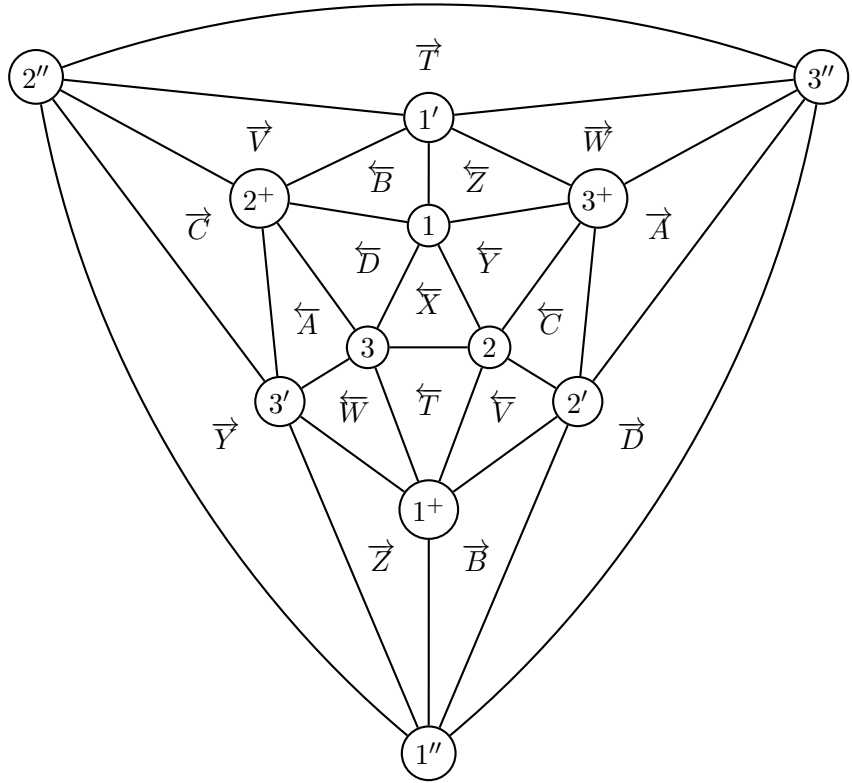
\centerline{Fig. 1: The icosahedron as a plane graph; the arrows indicate rotations}

\vspace{0.4cm}

In what follows, we will refer to Fig. 1 which represents an icosahedron drawn as a plane graph (for Graph Theory see \cite{H,Pl,sc}). We recall the presentation or defining relations of ${\cal A}_{5}$. Taking $D=(1,4,5)(2)(3)$, $Y=(2,4,5)(1)(3)$ and $T=(3,4,5)(1)(2)$ as generators we have $D^{3}=Y^{3}=T^{3}=I$ and $(DT)^{2}=(DY)^{2}=(YT)^{2}=I$ where $I$ is the identity element.

We now choose 3 rotations of the icosahedron that will also be designated by $D$, $Y$ and $T$, respectively (for the moment, forget the arrows!). In Fig. 1, they are represented by these same letters inside faces and with arrows indicating the direction of a $2\pi /3$ rotation. For example, $D$ is a rotation of 120 degrees around an axis joining the centers of the faces $\{1,2^{+},3\}$ and $\{1'',3'',2'\}$; it carries $1$ to $2^{+}$, $2^{+}$ to $3$, $3$ to $1$; and $1''$ to $2'$, $2'$ to $3''$, $3''$ to $1''$.

\vspace{0.2cm}

These rotations of 120 degrees around axes joining the central points of opposite faces are listed in Table 1. We call them {\it face rotations}. Remember that the identity may be regarded as a rotation of 0 (or 360) degrees!
$$\begin{array}{rccc cccc cccc cccc c}\
&1&2&3&&1'&2'&3'&&1^{+}&2^{+}&3^{+}&&1''&2''&3''\\\hline
D:&3&3'&2^{+}&&2&1''&1'&&2''&1&1^{+}&&3''&3^{+}&2'\\
Y:&3^{+}&1&1'&&2'&3&2''&&2^{+}&3''&2&&3'&1''&1^{+}\\
T:&2'&1^{+}&2&&3''&3'&1&& 3&3^{+}&1''&&2^{+}&1'& 2''\\
A:&1^{+}&1''&3'&&2&3''&2^{+}&&2''&3&2'&&1'&1&3^{+}\\
Z:&1'&2^{+}&2'' &&3^{+}&3&1''&&3'&3''&1&&1^{+}&2'&2\\
V:&3''&2'&3^{+}&&2''&1^{+}&1&&2&1'&1''&&3&2^{+}&3'\\
X:&2&3&1&&2'&3'&1'&&2^{+}&3^{+}&1^{+}&&2''&3''&1''\\
W:&2'&1''&1^{+}&&3^{+}&2''&3&&3'&2&3''&&2^{+}&1&1'\\
B:&2^{+}&3'&2''&&1&1^{+}&3''&&1''&1'&3&&2'&3^{+}&2\\
C:&3''&3^{+}&1'&&1''&2&2^{+}&&1&2''&2'&&3&3'&1^{+}
\end{array}$$

\centerline {Table 1: Face rotations and vertex permutations}

\vspace{0.2cm}

It is easy to verify that $D^{3}=Y^{3}=T^{3}=I$ and
$(DT)^{2}=(DY)^{2}=(YT)^{2}=I$; it is a bit more tedious to check that
$A=DYDT$, $Z=YTYD$, $V=TDTY$, $W=DTYT$, $B=YDTD$, $C=TYDY$ and
$X=YA^{2}YV$ (these products act left to right!).

To each one of these rotations, we associate a permutation of ${\cal A}_{5}$ of class $k_{1}=2,k_{3}=1$. It decomposes into one cycle of length 3 and two (trivial) cycles of length 1. We represent them just by the cycle of length 3 (We do not specify trivial cycles). Look that $(\alpha ,\gamma ,\beta)$ and
$(\alpha ,\beta ,\gamma )$ are inverse permutations. In the icosahedron they correspond to rotations of the same angle but in opposite directions. We have associated the icosahedron rotations $D,Y,T$ with the permutations $(1,4,5)$, $(2,4,5)$ and $(3,4,5)$, respectively. Now, see:
$$\begin{array}{lll}
(1,4,5)\leftrightarrow D&\ \ \ \ &(1,2,5)\leftrightarrow W\\
(2,4,5)\leftrightarrow Y&\ \ \ \ &(1,3,4)\leftrightarrow Z^{2}\\
(3,4,5)\leftrightarrow T&\ \ \ \ &(2,3,4)\leftrightarrow A\\
(1,3,5)\leftrightarrow C^{2}&\ \ \ \ &(1,2,4)\leftrightarrow V\\
(2,3,5)\leftrightarrow B&\ \ \ \ &(1,2,3)\leftrightarrow X\\
\end{array}$$

\centerline {Table 2: The 3-cycles of ${\cal A}_{5}$ and the
face rotations of the icosahedron}

\vspace{0.1cm}

If we recall how we wrote $A,Z,V,W,B,C,X$ as compositions of
$D,Y,T$, we recognize that Table 2 exhibits the bijection between permutations of class $k_{1}=2,k_{3}=1$ of ${\cal A}_{5}$ and face rotations of the icosahedron.

\vspace{0.2cm}

The next point is more delicate. We have to see how face rotations generate rotations of an angle $\pi $ around axes joining the middle points of opposite edges of the icosahedron: we call them {\it edge rotations}. We refer to Fig. 2.

Let $R_{1}$ be a face rotation that replaces vertices $x,y,z$ of a given face by the vertices $y,z,x$, respectively; let $R_{2}$ be another face rotation that replaces vertices $x,u,v$ of that other face by the vertices $u,v,x$, respectively. These faces share vertex $x$ and we can see them better if we look at the subgraph of the icosahedron induced by $x$ and its 5 neighbors. Such a subgraph is a wheel with $x$ as central vertex and $y,z,u,v,a$ as vertices of the peripheral pentagon.

\vspace{0.4cm}

\hspace{2.1cm}\psset{xunit=0.2cm,yunit=0.3cm,arrowlength=3,arrowsize=2pt
3} \cnodeput(0,0){y}{y}\cnodeput(10,0){z}{z}\cnodeput(30,0){z1}{z}
\cnodeput(40,0){y1}{y}
\cnodeput(0,-10){v}{v}\cnodeput(10,-10){u}{u}\cnodeput(30,-10){u1}{u}
\cnodeput(40,-10){v1}{v}
\cnodeput(-5,-5){a}{a}\cnodeput(5,-5){x}{x}\cnodeput(25,-5){a1}{a}
\cnodeput(35,-5){x1}{x}
\ncline{y}{z}\ncline{y}{x}\ncline{z}{x}\ncline{x}{v}
\ncline{x}{u}\ncline{u}{v}\ncline{a}{x}\ncline[doubleline=true]{z}{u}
\rput(5,-2){$\overleftarrow{R_{1}}$}\rput(5,-8){$\overleftarrow{R_{2}}$}
\ncline{y1}{z1}\ncline{y1}{x1}\ncline{z1}{x1}\ncline{x1}{v1}
\ncline{x1}{u1}\ncline{u1}{v1}\ncline[doubleline=true]{a1}{x1}\ncline{y1}{v1}
\rput(35,-2){$\overrightarrow{R_{1}}$}\rput(35,-8){$\overrightarrow{R_{2}}$}

\vspace{3.8cm}

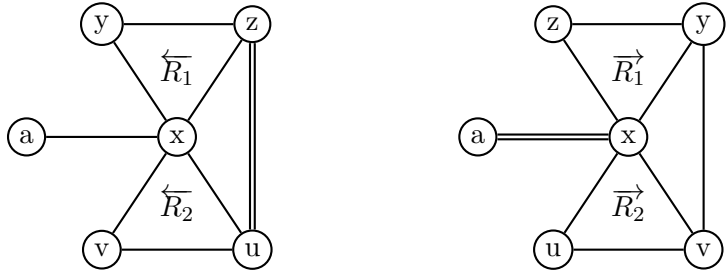
\centerline {Fig. 2: Face rotations generate edge rotations}

\vspace{0.4cm}

Suppose that, while moving along the pentagon, $y,z,u,v$ show up in this order. We distinguish two cases: either  $y,z,u,v,a$ or $y,z,a,u,v$. In the first case, the edge $[u,z]$ exists and the product  $R_{1}R_{2}$ is an edge rotation around the axis joining the middle point of $[u,z]$ to the middle point of the opposite edge; to make it short, we call it the edge $[u,z]$ rotation.

In the second case, the edge $[u,z]$ does not exist but there is an edge $[v,y]$ and $R_{1}R_{2}$ is an edge rotation around an axis joining the middle point of $[x,a]$ to the middle point of the opposite edge; $[x,a]$ is the only edge incident to $x$ and not belonging to either face $\{x,y,z\}$ or $\{x,u,v\}$. In Fig. 2, $[u,z]$ and $[x,a]$ are drawn in double line.

Note the following: if, while moving around the pentagon, vertices $y,z,v,u$ show up in this order, this means that $R_{1}$ and $R_{2}$ move the pentagon vertices in opposite directions along the pentagon. If this is so, then neither  $[u,z]$ nor $[v,y]$ exist and the product $R_{1}R_{2}$ is not an edge rotation.

We designate each edge rotation by the vertices of the respective pair of opposite edges. Table 3 associates to each edge rotation one ordered pair of face rotations that generate it. We designate the face rotations by the letters
they were given in Fig. 1. We will not explicitly give the whole vertex permutation which decomposes into 6 cycles of length 2; for instance, $CZ=(1,2)(1'',2'')(3,3^{+})(1',1^{+})(2',2^{+})(3',3'')$ but we indicate only the first two, in the example, $(1,2)(1'',2'')$ which are the vertices of the two opposite edges whose middle points are joined by the rotation axis.

Keep in mind that an edge rotation may be obtained by distinct pairs of face rotations. In the bottom line of Table 3, each one of the 3 rotations is generated by a configuration as in the right side of Fig. 2; the other 12 edge rotations of this table are generated by configurations as in the left side of Fig. 2. But each one of the 15 edge rotations may be obtained by products of face rotations as in the right or in the left side of Fig. 2.
$$\begin{array}{ccc}
(1,2),(1'',2'')= CZ&\ \ \ \ &(3^{+},2),(3',2'')= ZX\\
(2,3),(2'',3'')= YD&\ \ \ \ &(1^{+},2),(1',2'')= XW\\
(3,1),(3'',1'')= TY&\ \ \ \ &(2^{+},3),(2',3'')= XB\\
(1,3^{+}),(3',1'')= XC&\ \ \ \ &(1',3^{+}),(3',1^{+})= BY\\
(3,1^{+})(1',3'')= VX&\ \ \ \ &(1',2^{+}),(2',1^{+})= CT\\
(1,2^{+}),(2',1'')= AX&\ \ \ \ &(2',3^{+}),(3',2^{+})=
YV\\
\hline (1,1'),(1^{+},1'')= DY&\ (2,2'),(2^{+},2'')= YT&\
(3,3'),(3^{+},3'')= TD
\end{array}$$

\centerline {Table 3: Edge rotations as products of face rotations}

\vspace{0.2cm}

Now recall Table 2, and set up the correspondence between the 15 edge rotations of the icosahedron and the 15 permutations of ${\cal A}_{5}$ of class $k_{1}=1,k_{2}=2$; recall also the well known fact that $(\alpha ,\beta ,\gamma )(\alpha,\beta ,\delta)\equiv (\alpha ,\delta )(\beta ,\gamma)$. We get Table 4.
$$\begin{array}{lll}
CZ&\leftrightarrow (1,5,3)(1,4,3)\equiv (3,1,5)(3,1,4)&\equiv
(3,4)(1,5)(2)\\
YD&\leftrightarrow (2,4,5)(1,4,5)\equiv (4,5,2)(4,5,1)&\equiv
(4,1)(5,2)(3)\\
TY&\leftrightarrow (3,4,5)(2,4,5)\equiv (4,5,3)(4,5,2)&\equiv
(4,2)(5,3)(1)\\
XC&\leftrightarrow (1,2,3)(1,5,3)\equiv (3,1,2)(3,1,5)&\equiv
(3,5)(1,2)(4)\\
VX&\leftrightarrow (1,2,4)(1,2,3)&\equiv
(1,3)(2,4)(5)\\
AX&\leftrightarrow (2,3,4)(1,2,3)\equiv (2,3,4)(2,3,1)&\equiv
(2,1)(3,4)(5)\\
ZX&\leftrightarrow (1,4,3)(1,2,3)\equiv (3,1,4)(3,1,2)&\equiv
(3,2)(1,4)(5)\\
XW&\leftrightarrow (1,2,3)(1,2,5)&\equiv
(1,5)(2,3)(4)\\
XB&\leftrightarrow (1,2,3)(2,3,5)\equiv (2,3,1)(2,3,5)&\equiv
(2,5)(1,3)(4)\\
BY&\leftrightarrow (2,3,5)(2,4,5)\equiv (5,2,3)(5,2,4)&\equiv
(5,4)(2,3)(1)\\
CT&\leftrightarrow (1,5,3)(3,4,5)\equiv (5,3,1)(5,3,4)&\equiv
(5,4)(3,1)(2)\\
YV&\leftrightarrow (2,4,5)(1,2,4)\equiv (2,4,5)(2,4,1)&\equiv
(2,1)(4,5)(3)\\
DY&\leftrightarrow (1,4,5)(2,4,5)\equiv (4,5,1)(4,5,2)&\equiv
(4,2)(5,1)(3)\\
YT&\leftrightarrow (2,4,5)(3,4,5)\equiv (4,5,2)(4,5,3)&\equiv
(4,3)(5,2)(1)\\
TD&\leftrightarrow (3,4,5)(1,4,5)\equiv (4,5,3)(4,5,1)&\equiv
(4,1)(5,3)(2)
\end{array}$$

\centerline {Table 4: Edge rotations and permutations of class
$k_{1}=1,k_{2}=2$ of ${\cal A}_{5}$}

\vspace{0.1cm}

All that remains to be done is to analyze the so called {\it vertex rotations}, this means the 4 rotations (distinct from the identity) around each one of the 6 axes joining pairs of opposite vertices. These vertex rotations correspond to the 24 permutations of class $k_{5}=1$ of ${\cal A}_{5}$. For each one of the 6 axes, it is enough to express one of the 4 rotations in terms of the generators. We will do it for the rotations of $2\pi /5$ which appear in Figure 1 counterclockwise, since the other ones are powers of these ones. Let us take these rotations for the axes joining the vertices $\{1',1^{+}\}$, $\{2',2^{+}\}$,
$\{3',3^{+}\}$, $\{1,1''\}$, $\{2,2''\}$ e $\{3,3''\}$. We get
$$\begin{array}{llcc}
X^{2}D&=(1')(1^{+})&(1,2^{+},2'',3'',3^{+})&(2,3,3',1'',2')\\
X^{2}Y&=(2')(2^{+})&(1,1',2'',3',3)&(2,3^{+},3'',1'',1^{+})\\
X^{2}T&=(3')(3^{+})&(1,2,2',3'',1')&(3,1^{+},1'',2'',2^{+})\\
W^{2}A&=(1)(1'')&(2,3,2^{+},1',3^{+})&(2',1^{+},3',2'',3'')\\
B^{2}Z&=(2)(2'')&(1,3^{+},2',1^{+},3)&(1',3'',1'',3',2^{+})\\
C^{2}V&=(3)(3'')&(1,2,1^{+},3',2^{+})&(1',3^{+},2',1'',2'')
\end{array}$$
\centerline {Table 5: Vertex rotations generated by face rotations}

\vspace{0.2cm}

Looking at Fig. 1, we see that, for instance, in the rotation $W^{2}A$, the vertices adjacent to the fixed vertex $1$ rotate around it counterclockwise; vertices $3,2^{+},1',3^{+}$ and $2$ are replaced by
vertices $2^{+},1',3^{+},2$ and $3$, respectively.

To each vertex rotation we now associate a permutation of class $k_{5}=1$ of ${\cal A}_{5}$.
Denote by $S_{i}$, for $i=1,...,6$, the rotations which appear in Table 5 and let us associate to each one of them a permutation $Q_{i}$ of ${\cal A}_{5}$ which will be determined based on the correspondences of Table 2. We get

$$\begin{array}{l}
S_{1}=X^{2}D\leftrightarrow (1,3,2)(1,4,5)=(1,3,2,4,5)=Q_{1}\\
S_{2}=X^{2}Y\leftrightarrow (1,3,2)(2,4,5)=(1,3,4,5,2)=Q_{2}\\
S_{3}=X^{2}T\leftrightarrow (1,3,2)(3,4,5)=(1,4,5,3,2)=Q_{3}\\
S_{4}=W^{2}A\leftrightarrow (1,5,2)(2,3,4)=(1,5,3,4,2)=Q_{4}\\
S_{5}=B^{2}Z\leftrightarrow (2,5,3)(1,4,3)=(1,4,3,2,5)=Q_{5}\\
S_{6}=C^{2}V\leftrightarrow (1,3,5)(1,2,4)=(1,3,5,2,4)=Q_{6}
\end{array}.$$

\centerline {Table 6: Permutations of ${\cal A}_{5}$ and $2\pi /5$ rotations of the icosahedron}

\vspace{0.2cm}

We end up with the correspondence between the remaining permutations of class $k_{5}=1$ of ${\cal A}_{5}$ (there are 18 besides the 6 already indicated) and the remaining vertex rotations in angles of $4\pi /5$, $6\pi /5$ and $8\pi /5$ counterclockwise.
$$\begin{array}{lll}
(1,2,3,4,5)=Q_{6}^{3}&\ (1,3,4,2,5)=Q_{3}^{3}&\
(1,4,5,2,3)=Q_{4}^{3}\\
(1,2,3,5,4)=Q_{3}^{4}&\ (1,3,4,5,2)=Q_{2}&\
(1,4,5,3,2)=Q_{3}\\
(1,2,4,3,5)=Q_{4}^{4}&\ (1,3,5,2,4)=Q_{6}&\
(1,5,2,3,4)=Q_{5}^{4}\\
(1,2,4,5,3)=Q_{5}^{3}&\ (1,3,5,4,2)=Q_{5}^{2}&\
(1,5,2,4,3)=Q_{3}^{2}\\
(1,2,5,3,4)=Q_{1}^{2}&\ (1,4,2,3,5)=Q_{2}^{2}&\
(1,5,3,2,4)=Q_{2}^{3}\\
(1,2,5,4,3)=Q_{2}^{4}&\ (1,4,2,5,3)=Q_{6}^{4}&\
(1,5,3,4,2)=Q_{4}\\
(1,3,2,4,5)=Q_{1}&\ (1,4,3,2,5)=Q_{5}&\
(1,5,4,2,3)=Q_{1}^{4}\\
(1,3,2,5,4)=Q_{4}^{2}&\ (1,4,3,5,2)=Q_{1}^{3}&\
(1,5,4,3,2)=Q_{6}^{2}
\end{array}$$

\centerline {Table 7: To the permutation $Q_{\alpha }^{\beta }$ of
${\cal A}_{5}$ we associate the rotation $S_{\alpha }^{\beta }$}

\vspace{0.2cm}

Our last point is trivial. It is immediate to check that the 24 permutations of ${\cal A}_{5}$ in this class are the powers of $Q_{i}$ with exponents 1, 2, 3 and 4. They correspond one-to-one to the powers of the same exponent of the rotations $S_{i}$. In Table 7, we present these 24 permutations in lexicographic order and how they can be written as powers of each $Q_{i}$.

\vspace{0.2cm}

(*) Department of Mathematics, University of Coimbra, Coimbra, Portugal

URL: www.mat.uc.pt

Author's e-mail: siper@mat.uc.pt

\end{document}